\documentclass[12pt,centertags,oneside]{amsart}
\usepackage{amsmath,amstext,amsthm,amscd,typearea,hyperref}
\usepackage{amssymb}
\usepackage{a4wide}
\usepackage[mathscr]{eucal}
\usepackage{mathrsfs}
\usepackage{typearea}
\usepackage{charter}
\usepackage{pdfsync}
\usepackage[a4paper,width=16.2cm,top=3cm,bottom=3cm]{geometry}

\numberwithin{equation}{section}

\newtheorem{theorem}{Theorem}[section]

\newtheorem{proposition}[theorem]{Proposition}
\newtheorem{corollary}[theorem]{Corollary}
\newtheorem{lemma}[theorem]{Lemma}
\newtheorem{remark}[theorem]{Remark}

\newcommand{\cali}[1]{\mathscr{#1}}

\DeclareMathOperator{\Tan}{Tan}
\DeclareMathOperator{\Per}{Per}
\DeclareMathOperator{\IPer}{IPer}

\newcommand{\dbar}{\overline\partial}

\newcommand{\Cc}{\cali{C}}

\newcommand{\Oc}{\cali{O}}

\newcommand{\C}{\mathbb{C}}

\newcommand{\N}{\mathbb{N}}
\newcommand{\Z}{\mathbb{Z}}
\newcommand{\R}{\mathbb{R}}
\renewcommand\P{\mathbb{P}}

\newcommand{\rI}{{\textrm{\rm I}}}
\newcommand{\rII}{{\textrm{\rm II}}}

\title[Growth of the number of   periodic points]{Growth of the number of   periodic points for meromorphic maps}

\author{Tien-Cuong Dinh}
\address{Department of Mathematics, National University 
of Singapore, 10 Lower Kent Ridge Road, Singapore 119076.}
\email{matdtc@nus.edu.sg}
\thanks{T.-C.\ D.\  was supported by Start-Up 
Grant R-146-000-204-133 from National University of
\break Singapore}
\author{Vi{\^e}t-Anh Nguy{\^e}n}
\address{Math{\'e}matique-B{\^a}timent 425, UMR 8628, 
Universit{\'e} Paris-Sud, 91405 Orsay, France.}
\email{VietAnh.Nguyen@math.u-psud.fr}
\author{Tuyen Trung Truong}
\address{Department of Mathematics, School of Mathematical Sciences, The University of Adelaide, Adelaide, SA 5005, Australia.}
\email{tuyen.truongt@adelaide.edu.au}
\thanks{T.T. T.\  was partially supported by  Australian Research Council grants DP120104110 and DP150103442.}

\date{January 12, 2016}
\begin{document}

\begin{abstract}
We show that any dominant meromorphic self-map $f:X\rightarrow X$ of a  compact K{\"a}hler
manifold $X$ is an Artin-Mazur map. More precisely, if  $P_n(f)$ is the number of  its isolated  periodic points of period $n$ (counted with multiplicity), then  
$P_n(f)$ grows at most exponentially fast with respect to $n$ and the exponential rate is at most equal to the 
algebraic entropy  of $f$.
Further estimates are given when $X$ is a  surface. Among the techniques introduced in this paper,  the h-dimension of the density between two arbitrary positive closed currents on a compact K\"ahler surface is obtained.
\end{abstract}

\maketitle

\medskip

\noindent
{\bf Classification AMS 2010}: 37F, 32U, 32H50.

\medskip


\section{Introduction} \label{introduction}

A self-map on a compact differentiable manifold is called {\it an Artin-Mazur map} if its number of isolated periodic points of period $n$ grows at most exponentially fast with $n$. Artin and Mazur proved in \cite{ArtinMazur} that such maps are dense in the set of $\Cc^k$ maps. On the other hand, Kaloshin constructed large families of diffeomorphisms such that the number  of isolated periodic points of period $n$ grows faster than any given sequence of integers \cite{Kaloshin}. In particular, there are maps which do not satisfy Bowen's formula relating the topological entropy  and the exponential growth rate of number of periodic points. Furthermore, the dynamical $\zeta$-function, associated with such a map,  is not analytic in any neighborhood of zero. 

In this paper, we consider the question in the complex setting.
Let $X$ be a  compact K\"ahler manifold of dimension $k$.
Let  $f:\ X\to X$ be a meromorphic self-map of $X$. 
We always assume  that $f$ is {\it dominant,}  i.e.,    the image of $f$ contains an open subset of $X$.  [Otherwise, the study of this map is reduced to the case of dominant maps in a lower dimensional manifold. 
Note also that meromorphic maps may have indeterminacy sets and hence are not continuous in general.] 
We will show that such a map is always Artin-Mazur and an upper 
bound for the growth rate of the number of periodic points, in the spirit of Bowen's formula, holds in general. More precise results will be given in the case of dimension $k=2$. 

In order to state our main results, let us recall some basic notions. The map $f$ is holomorphic outside a (possibly empty) analytic set of co-dimension at least 2 which is called {\it the indeterminacy set} of $f$ and is denoted by $I(f)$. The graph of $f$ over $X\setminus I(f)$ can be compactified to be an irreducible analytic set of dimension $k$ in $X\times X$ and we denote this analytic set by $\Gamma(f)$. 
If $\pi_1,\pi_2$ are the natural projections from $X\times X$ to its factors and $[\Gamma(f)]$ the current of integration on $\Gamma(f)$, we define the pull-back operator $f^*$ on smooth differential $(p,q)$-forms $\phi$ by 
$$f^*(\phi):=(\pi_1)_*\big(\pi_2^*(\phi)\wedge [\Gamma(f)]\big).$$
This is a $(p,q)$-current defined by an $L^1$ form which is smooth outside $I(f)$. The operator commutes with $\partial,\dbar$ and therefore induces a natural pull-back operator $f^*:H^{p,q}(X,\C)\to H^{p,q}(X,\C)$ on the Hodge cohomology group $H^{p,q}(X,\C)$ of $X$.

The iterate of order $n$ of $f$ is defined by $f^n:=f\circ\cdots\circ f$
($n$ times) on a dense Zariski open set and extends to a
dominant meromorphic map on $X$.  Since the above current $f^*(\phi)$  is not smooth in general, we cannot iterate the operators on differential forms. One can consider the iterate of $f^*$ acting on Hodge cohomology but we don't have the identity $(f^n)^*=(f^*)^n$  for every $f$. This well-known phenomenon, observed by Forn\ae ss and Sibony, is a source of difficulties in the dynamical study of $f$.  However, for $0\leq p\leq k$ and for any fixed norm on cohomology, the norm of $(f^n)^*$ on $H^{p,p}(X,\C)$ has a nice behavior : it was shown by Sibony and the first author that the limit 
$$d_p(f):=\lim_{n\to\infty} \|(f^n)^*\|_{H^{p,p}(X,\C)}^{1/n}$$
always exists and is a fundamental bi-meromorphic (finite) invariant of $f$, independent of the choice of the norm on cohomology. This is {\it the dynamical degree of order $p$} of $f$. We always have $d_0(f)=1$, $d_p(f)\geq 1$, $d_p(f^n)=d_p(f)^n$, and that $p\mapsto \log d_p(f)$ is concave, i.e., $d_p(f)^2\geq d_{p-1}(f)d_{p+1}(f)$ for $1\leq p\leq k-1$.  The last dynamical degree $d_k(f)$ is also called {\it the topological degree}. It is equal to the number of points in a generic fiber of $f$. 
{\it The algebraic entropy} of $f$ is defined by
$$h_a(f):=\max_{0\leq p\leq k} \log d_p(f).$$
The topological entropy of $f$ is always bounded above by the algebraic entropy, see \cite{DinhNguyenTruong1, DinhSibony04, DinhSibony05, Gromov, Yomdin} for details.

Let $\Delta$ denote the diagonal of $X\times X$ and $\Gamma(f^n)$ the closure of the graph of $f^n$ in $X\times X$. A point $x\in X$ 
is a {\it periodic point of period $n$} of $f$ if $(x,x)\in \Delta \cap \Gamma(f^n)$, and such a point is {\it isolated} if $(x,x)$ is an isolated point in $\Delta \cap \Gamma(f^n)$.
 Note that the set of periodic points of period $n$ of $f$ may have positive dimension even when $f$ is a holomorphic map.
 Denote by $\Per_n(f)$ (resp. $\IPer_n(f)$)  the set of all periodic (resp. isolated periodic)  points 
  of period $n$ of $f.$  
For $x\in \IPer_n(f),$ the local index of $f^n$ at $x,$ denoted by  $\nu_x(f^n),$ is   the local multiplicity
of the intersection $ \Delta \cap \Gamma(f^n)$ at $(x,x)$, see e.g. \cite{Mumford}.
Let
$$P_n(f):=\sum_{x\in \IPer_n(f)}  \nu_x(f^n)$$
be the number of isolated periodic points of period $n$ counted with multiplicity.
The following dynamical $\zeta$-function of $f$ is similar to the one introduced by Artin and Mazur in \cite{ArtinMazur}
$$\zeta_f(z):=\exp\Big(\sum_{n\geq 1} P_n(f){z^n\over n}\Big).$$

Here is our first main theorem.

\begin{theorem}\label{th_main_1}
Let $f$ be a dominant  meromorphic self-map  on a compact K\"ahler manifold $X$.  Let $h_a(f)$ be its algebraic entropy and $P_n(f)$ its number of isolated periodic points of period $n$ counted with multiplicity. 
Then we have 
$$\limsup_{n\to\infty}{1\over n} \log  P_n(f) \leq h_a(f).$$
In particular, $f$ is an Artin-Mazur map, i.e., its number of  isolated periodic points of period $n$ grows at most exponentially fast  with $n$, and therefore, the dynamical $\zeta$-function of $f$ is analytic in a neighborhood of zero. 
 \end{theorem}

Note that the above result can be extended to meromorphic correspondences with essentially the same proof, see also \cite{DinhSibony08}.
When $f$ has only isolated periodic points, i.e., $\IPer_n(f)=\Per_n(f)$, the number $P_n(f)$ is equal to the intersection number between $\Gamma(f^n)$ and $\Delta$ which can be computed, via the Lefschetz fixed point formula, in terms of the trace of the operator $(f^n)^*$ on Hodge cohomology. In this case, the above theorem can be easily obtained. Also, when $X$ is homogeneous, one can move slightly $\Delta$ so that its intersection with $\Gamma(f^n)$ is finite and see that the number of isolated intersection points cannot decrease. An upper bound for $P_n(f)$ can be easily obtained in this case. 

Note that if $\pi:X'\to X$ is a bi-meromorphic map between compact K\"ahler manifolds, then we can lift $f$ to a dominant meromorphic self-map 
of $X'$ by setting $f':=\pi^{-1}\circ f\circ \pi$. The maps $f$ and $f'$ induce very related dynamical systems but we don't have $P_n(f)=P_n(f')$ in general. In other words, $P_n(f)$ is not a bi-meromorphic invariant of $f$. So even when $X$ is bi-meromorphic to a homogeneous manifold, one cannot directly reduce the problem to this case.

In order to obtain the theorem in the general case, we will use a recent theory of densities of currents introduced by Sibony and the first author \cite{DinhSibony12}. Roughly speaking, this theory allows us to dilate the coordinates in the normal directions to $\Delta$ in $X\times X$ by a factor $\lambda$. Taking $\lambda\to\infty$, the limit of the image of $\Gamma(f^n)$ by the dilation allows us to bound the number of isolated periodic points in terms of the volume of $\Gamma(f^n)$. The last quantity can be estimated using dynamical degrees of $f$.

We also expect, at least for large families of maps, that the numbers $P_n(f)$ satisfy
$$P_n(f)=e^{nh_a(f)}+o(e^{nh_a(f)}) \quad \text{or equivalently} \quad \lim_{n\to\infty} e^{-nh_a(f)} P_n(f)=1.$$
The problem is widely open. Different approaches were introduced to construct a good number of repelling or saddle periodic points for particular families of maps. This allows for obtaining a good lower bound for $P_n(f)$. We refer to \cite{BedfordLyubichSmillie, BriendDuval99, DillerDujardinGuedj3, DinhNguyenTruong,DinhSibony10, DinhSibony14, Dujardin} for results in this direction.

We will now concentrate in the case where $X$ is a compact K\"ahler surface. Despite deep knowledges about complex surfaces, the above basic question is still open. We refer to \cite{Favre,IwasakiUehara, Saito, Xie} for some results related to this question. 
When the topological degree of $f$ is dominant, i.e., $d_2(f)>d_1(f)$, the affirmative answer was recently obtained by the authors in \cite{DinhNguyenTruong}. From now on, we assume that $f$ has {\it minor topological degree}\footnote{Such a map is called a map with small topological degree in literature. We think that our terminology is more appropriate.} in the sense that $d_2(f)<d_1(f)$. In this case, we have $h_a(f)=\log d_1(f)$. 
We also assume that $f$ is {\it algebraically stable} in the sense of Forn\ae ss-Sibony, i.e., $(f^n)^*=(f^*)^n$ as linear maps on $H^{1,1}(X,\mathbb{C})$ for all $n\in \mathbb{N}$. This condition can be checked geometrically. 
Here is our second main result. 

\begin{theorem}\label{th_main_2}
Let $f$ be a meromorphic self-map  with  minor topological degree on a compact K\"ahler surface $X$. Let $d_1(f)$ denote the first dynamical degree and $P_n(f)$ the number of isolated periodic points of period $n$ of $f$ counted with multiplicity. Assume that  $f$ is  algebraically stable in the sense of Forn\ae ss-Sibony. Then we have
$$P_n(f)\leq d_1(f)^n + o( d_1(f)^n ) \quad \text{as} \quad n\to\infty.$$
\end{theorem}

Note that in this case, the Lefschetz number, i.e., the intersection number between $\Gamma(f^n)$ and $\Delta$, is not difficult to estimate. It is also equal to $ d_1(f)^n + o( d_1(f)^n )$, see Lemma \ref{L:Lefschetz} below.  Observe that saddle periodic points are isolated and of multiplicity 1. The following result is an immediate consequence of the last theorem.

\begin{corollary} \label{cor_main}
Under the conditions of Theorem \ref{th_main_2}, assume moreover that 
there is a family $A_n$ of saddle periodic points of period $n$ of $f$ of cardinality $d_1(f)^n+o(d_1(f)^n)$ which are equidistributed with respect to a probability measure $\mu$, i.e., 
$$\lim_{n\to \infty}  d_1(f)^{-n}\sum_{a\in A_n} \delta_a = \mu,$$
where $\delta_a$ stands for the Dirac mass at $a$. Then the set of all saddle periodic points and the set of all isolated periodic points (counting multiplicity or not) of period $n$ are also of cardinality  $d_1(f)^n+o(d_1(f)^n)$ and equidistributed with respect to $\mu$ as $n\to\infty$. 
\end{corollary}

Under technical conditions, by using Pesin's theory, saddle periodic points can be constructed and satisfy the hypothesis of the last corollary. We refer to \cite{DillerDujardinGuedj3, Dujardin, JonssonReschke} for more precision on these conditions, see also \cite[pp.694-5]{KatokHasselblatt}. Corollary \ref{cor_main} emphasizes 
the role of the upper bound of the number of isolated periodic points (Theorem \ref{th_main_2}) in their equidistribution property. This upper bound is sometimes overlooked in literature.

Our last main result  deals with algebraically stable bi-meromorphic surface maps $f:\ X\to X$. 
To this end, we will make use of the    Saito's local index function  $\nu_\bullet(f^n): \Per_n(f)\to\N,$ 
introduced  by  Saito  \cite{Saito} and  Iwasaki-Uehara \cite{IwasakiUehara},
 which extends the usual local index function on isolated periodic points to all periodic points, see Section  \ref{S:bimero_maps}  below for more details.
 Define   
 $$P'_n(f):=\sum_{x\in \Per_n(f)}  \nu_x (f^n).$$  
 Clearly, $P_n(f) \leq P'_n(f).$

\begin{theorem}\label{th_main_3}
Let $f$ be a bi-meromorphic self-map on a compact K\"ahler surface $X$ which is is  algebraically stable. Assume that its first dynamical degree satisfies $d_1(f)>1$.  Then 
$$P'_n(f)= d_1(f)^n + o( d_1(f)^n ) \quad \text{as} \quad n\to\infty.$$
\end{theorem}

Note that in dynamics one often assumes that $d_1(f)>1$ because otherwise the topological entropy of $f$ is zero and, in some sense, the dynamics is poor. 

The paper is  organized  as  follows.
In Section  \ref{S:currents} we  recall some  background on positive  closed currents and the theory of densities for currents. 
Basic properties of the action of maps on currents and cohomology will be presented in 
 Section \ref{S:growth_rate} together with the proofs of Theorems \ref{th_main_1} and  \ref{th_main_2}. 
A main ingredient is the theory of densities for currents. Knowledge on dynamical Green currents of $f$ is used to get sharp upper bounds on $P_n(f)$ in the second theorem.
Finally, the case of bi-meromorphic maps together with Saito's and Iwasaki-Uehara's index
will be presented in  Section \ref{S:bimero_maps}.
This  key point in the proof of Theorem \ref{th_main_3}  is the counterpart of Lefschetz fixed point formula in the presence of non-isolated periodic points. Using  the theory of tangent currents as well as  some  recent results on periodic curves, we are able 
to control   the contribution of non-isolated periodic points.

\medskip
\noindent
{\bf Acknowledgement.} The paper was partially written during the visits of the second author at National University of Singapore and  Vietnam Institute for Advanced Study in Mathematics (VIASM). He would like to thank these organizations for
hospitality and support. The third author would like to thank Nicholas Buchdahl,  Charles Favre and Takano Uehara for fruitful discussions on this research topic.


\section{Positive closed currents and theory of densities}       \label{S:currents}

In this section, we first recall some basic properties of  positive closed currents on a  compact K\"ahler manifold. We then present some properties of tangent currents, a fundamental notion in theory of densities, that will be used in this work. The reader will find more details in  \cite{Demailly, DinhSibony10, DinhSibony12,Voisin}.
 
Let $X$ be a compact K\"ahler manifold of dimension $k$ and $\omega$ a fixed K\"ahler form on $X$.
Denote by 
$H^{p,q}(X,\C)$ the Hodge cohomology group of bi-degree $(p,q)$ of $X$ and define $H^{p,p}(X,\R):=H^{p,p}(X,\C)\cap H^{2p}(X,\R)$.
The cup-product in $\oplus H^*(X,\C)$ is  denoted by $\smallsmile.$ 
If $T$ is a closed $(p,q)$-current, its class in $H^{p,q}(X,\C)$ is denoted by $\{T\}$. If $T$ is a positive closed $(p,p)$-current, its class belongs to $H^{p,p}(X,\R)$. 
Recall  that  the pseudo-effective cone $H^{p,p}_{psef}(X)\subset H^{p,p}(X,\R)$ is the set of cohomology classes of positive closed
$(p,p)$ currents. It is  closed, convex and salient in the sense that it contains no non-trivial vector subspace. 
The nef cone $H^{1,1}_{nef}(X)$ is the closure of the cone of  all K\"ahler classes, i.e., classes of K\"ahler forms. It is contained in the pseudo-effective cone $H^{1,1}_{psef}(X)$. So it is also convex, closed and salient.

If $T$ is a current on $X$ and $\varphi$ is a test form of complementary degrees, the pairing $\langle T,\varphi\rangle$ denotes the value of $T$ at $\varphi$.  If $T$ is a positive $(p,p)$-current on $X$, we use the following notion of mass for $T$  given by the formula
$$\|T\|:=\langle T,\omega^{k-p}\rangle.$$
This mass is equivalent to the usual mass for order 0 currents. 
When $T$ is positive and closed, its mass depends only on its cohomology class $\{T\}$ in $H^{p,p}(X,\R)$. This is a key point in the calculus with positive closed currents. 
We will write $T\leq T'$ and $T'\geq T$ for two real $(p,p)$-currents $T,T'$ if $T'-T$ is a positive current. We also write $c\leq c'$ and $c'\geq c$ for $c,c'\in H^{p,p}(X,\R)$ when $c'-c$ is the class of a positive closed $(p,p)$-current, i.e., an element of $H^{p,p}_{psef}(X,\R)$. If $V$ is an analytic subset of pure dimension $k-p$ in $X$, denote by $[V]$ the positive closed $(p,p)$-current of integration on $V$ and $\{V\}$ its cohomology class in $H^{p,p}(X,\R)$. 
 Given a  positive closed current $T$ and a point $x\in X,$
  the Lelong number of $T$ at $x$ is denoted by $\nu(T,x)$.

We recall now basic facts about tangent currents and prove some abstract results which will allow us to bypass Lefschetz fixed point formula in order to bound the number of periodic points. We will restrict ourselves to the simplest situation that is  needed for the present work, see \cite{DinhSibony12} for more details.

Let $V$ be an irreducible  submanifold of $X$ of  dimension $l$ and
$\pi:\ E\to V$ be the normal vector bundle  of $V$ in $X$.
For  a point $a\in V,$ if   $\Tan_a X$ and $\Tan_a V$ denote respectively the tangent spaces of $X$ and $V$ at $a,$
the fiber $E_a:=\pi^{-1}(a)$  of $E$ over $a$ is canonically identified with the quotient space 
$\Tan_a X / \Tan_a V.$ The zero section of $E$  is naturally identified with $V$. 
Denote by $\overline E$ the natural compactification of $E$, i.e.,
the projectivization $\P(E\oplus \C)$ of the vector bundle $E\oplus \C$, where $\C$ is the trivial line bundle over $V$.
We still denote by $\pi$ the natural projection from $\overline E$ to $V$. 
Denote by $A_\lambda$ the multiplication by $\lambda$ on the fibers of  $E$ for $\lambda\in\C^*,$
i.e., $A_\lambda(u):=\lambda u$  for $u\in E_a$ and  $a\in V.$ This map extends to a holomorphic automorphism of $\overline E$. 

Let $V_0$ be an open subset of $V$ which is naturally  identified with an open subset of the section 0 in $E$. A diffeomorphism $\tau$ from a neighborhood of $V_0$ in $X$ to a neighborhood of $V_0$ in $E$ is called 
{\it admissible} if it satisfies  essentially the following three conditions:
the restriction of $\tau$ to $V_0$ is the  identity,  the differential  of $\tau$ at each point  $a\in V_0$  is
$\C$-linear and the composition  of the following maps
$$E_a\hookrightarrow \Tan_a(E)\to \Tan_a(X)\to E_a$$    
is the identity map on $E_a$. Here, the morphism $\Tan_a(E)\to \Tan_a(X)$ is given by the differential  of $\tau^{-1}$ at   $a$ and the other maps are the canonical ones, see  \cite{DinhSibony12}.

When $V_0$ is small enough, 
 there are local holomorphic coordinates on a small neighborhood $U$ of $V_0$ in $X$ so that over $V_0$ we identify naturally $E$ with $V_0\times \C^{k-l}$ and $U$ with an open neighborhood of $V_0\times\{0\}$  in $V_0\times \C^{k-l}$ (we reduce $U$ if necessary). 
 In this picture, the identity is a holomorphic admissible map.
 This  picture  is  called a {\it standard local setting}. 
Note  that a general admissible map is not necessarily holomorphic and
 there always  exist admissible maps for $V_0:=V$, which are rarely  holomorphic.  

Consider an admissible map $\tau$ as above.
Let $T$ be a positive closed $(p,p)$-current on $X$ without mass on $V$ for simplicity.
Define 
$$T_\lambda:=(A_\lambda)_*\tau_*(T).$$
The family $(T_\lambda)$ is relatively compact on $\pi^{-1}(V_0)$ when $\lambda\to\infty$: we can extract convergent subsequences for $\lambda\to\infty$. The limit currents $R$  are positive closed $(p,p)$-currents on $\overline E$ without mass on $V$. They are $V$-conic, i.e.,
$(A_\lambda)_* R =R$ for any $\lambda\in\C^*,$ or equivalently, $R$ is 
invariant by $A_\lambda.$  

Such a current $R$ depends on the choice of a subsequence of  $(T_\lambda)$ but it is independent of the choice of $\tau$. This property gives us a large flexibility to work with admissible maps.
In particular, using global admissible maps, we obtain positive closed $(p,p)$-currents $R$ on $\overline E$ though $\tau$ is not holomorphic. 
It is also known that the cohomology class of $R$ depends on $T$ but does not depend on the choice of $R$. 
This class is denoted by $\kappa^V(T)$ and called {\it the total tangent class} of $T$ with respect to $V$. The currents $R$ are {\it the tangent currents} of $T$ along $V$.  The mass of $R$ and the norm of $\kappa^V(T)$ are bounded by a constant times the mass of $T$.

Let $-h$ denote the tautological $(1,1)$-class on $\overline E$. Recall that $\oplus H^*(\overline E,\C)$ is a free $\oplus H^*(V,\C)$-module generated by $1,h,\ldots, h^{k-l}$ (the fibers of $\overline E$ are of dimension $k-l$). So we can write {\bf in a unique way}

\begin{equation}\label{e:total_class}
\kappa^V(T)=\sum\limits_{j=\max(0,l-p)}^{\min(l,k-p)} \pi^*(\kappa^V_j(T))\smallsmile h^{j-l+p}
,\end{equation}
where  $\kappa^V_j(T)$ is  a  class in $H^{l-j,l-j}(V,\R)$ with the convention  that $\kappa^V_j(T)=0$  outside the range
$ \max(0,l-p)\leq j\leq  \min(l,k-p) .$

The maximal integer $s$ such that $\kappa_s^V(T)\not=0$ is called {\it the tangential h-dimension} of $T$ along $V$
(when the tangent currents of $T$ along $V$ vanish, such an integer doesn't exist and we set the tangential h-dimension equal to 0).  
The class $\kappa^V_s(T)$ is pseudo-effective, i.e., contains a positive closed $(l-s,l-s)$-current on $V.$ 
Since we assume that $T$ has no mass on $V$, we always have $s<k-p$. 
 If $\omega_V$ is any K\"ahler form on $V$, 
the  tangential h-dimension of $T$ is  also  equal to the maximal integer $s\geq 0$ such that $R\wedge \pi^*(\omega_V^s)\not=0$, except when $R=0$.
In particular,  when $V$ is reduced to a single point $x,$ we have $s=0$ and $\kappa^{x}_0(T)=\nu(T,x)\{x\}.$
So  the notion of total tangent classes  generalizes
the notion of Lelong number. The next result may be  regarded as  the counterpart of Siu's semi-continuity theorem for Lelong numbers. 
 
 \begin{proposition}\label{P:semicontinuity} {\rm  (\cite[Theorem 4.11]{DinhSibony12})}
With the above notation, let  $T_n$ and $T$ be positive closed $(p,p)$-currents on $X$ such that
 $T_n\to T$. Let $s$ be  an integer at least equal to the  tangential h-dimension of $T$ along $V$. Then $\kappa_j^V(T_n)\to 0$ for $j>s$ and any limit class of $\kappa_s^V(T_n)$ is  pseudo-effective and is smaller than or equal to $\kappa_s^V(T)$. 
\end{proposition}

The following result will allow us to bound the number of isolated periodic points of a meromorphic map. We always identify the cohomology group $H^{2k}(X,\C)$ with $\C$ using the integrals of top degree differential forms on $X$. 
 
 \begin{proposition} \label{prop_intersection}
 Let $\Gamma_n$ be complex subvarieties of pure dimension $k-l$ in $X$. Assume that there is a sequence of positive numbers $d_n$ such that $d_n\to \infty$ and $d_n^{-1}[\Gamma_n]$
 converges to a positive closed $(l,l)$-current $T$ on $X$. Assume also that the tangential  h-dimension of $T$ with respect to $V$ is $0$ and that $\{T\}\smallsmile \{V\}=c\in\R^+$. Then the number 
$\delta_n$ of isolated points in the intersection $\Gamma_n\cap V$, counted with multiplicity, satisfies $\delta_n\leq cd_n+o(d_n)$ as $n\to\infty$. 
 \end{proposition}
 \proof 
 The case  $c=1$ has been proved in \cite[Proposition 2.1]{DinhNguyenTruong} using the above semi-continuity for densities.
 The proof there works for all $c\geq 0.$ We can also deduce the general case directly from the case $c=1$. 
 \endproof

 Let $T$ and $S$ be positive closed currents on $X$. The density between $T$ and $S$ is represented by the tangent currents of $T\otimes S$ along the diagonal $\Delta$ of $X\times X$. 
We will not give here the formal definition of this density and refer to \cite{DinhSibony12} for details. What we will need is the tangential  h-dimension of $T\otimes S$ with respect to $\Delta$. We call it {\it the h-dimension of the density between $T$ and $S$.}
In the special case when $S$ is a positive measure, the density between $T$ and $S$ can be represented by a number which is the integral of the Lelong number $\nu(T,x)$ with respect to the measure $S$. In this case, the density is always of h-dimension 0 and doesn't vanish if and only if  the above integral of $\nu(T,x)$ with respect to $S$ is non-zero. 
We have the following proposition.

 \begin{proposition} \label{P:tensor_prod_hdim}
Let $T$ and $S$ be two positive closed currents on $X$ of bi-degrees $(p,p)$ and $(q,q)$ respectively with $1\leq p,q\leq k-1$.
If $S$ has no mass on the set $\{x\in X,\ \nu(T,x)>0\}$,   
then the h-dimension of the density between $T$ and $S$ is at most equal to  $k-q-1$. 
 \end{proposition} 
\proof 
Consider the positive measure $S':=S\wedge \omega^{k-q}$. By hypothesis, $S'$ has no mass on $\{x\in X,\ \nu(T,x)>0\}$. 
We deduce that the density between $T$ and $S'$ vanishes. 
Assume by contradiction that the h-dimension of the density between $T$ and $S$ is larger than or equal to $k-q$. 
We will obtain a contradiction by showing that the density between $T$ and $S'$ doesn't vanish. 

Let $\pi_1$ and $\pi_2$ denote the natural projections from $X\times X$ to its factors. Similarly to the beginning of this section, denote  by $\pi:E\to\Delta$ the normal vector bundle to $\Delta$ in $X\times X$ which extends to the natural compactification $\overline E$ of $E$. Let $A_\lambda:\overline E\to\overline E$ be the multiplication by $\lambda$ on the fibers of $\pi$. Fix also a sequence $\lambda_n\to\infty$ so that 
$(A_{\lambda_n})_*\tau_*(T\otimes S)$ converges to a tangent current $R$ of $T\otimes S$ along $\Delta$ which does not depend on the choice of any admissible map $\tau$ from a neighborhood of $\Delta$ in $X\times X$ to a neighborhood of $\Delta$ in $E$.  

If $\alpha$ is any continuous positive closed form on $X\times X$, it is easy to deduce that $(A_{\lambda_n})_*\tau_*\big((T\otimes S)\wedge\alpha\big)$ converges to $R\wedge \pi^*(\alpha_{|\Delta})$. In particular, we have 
$$(A_{\lambda_n})_*\tau_*(T\otimes S') \to R\wedge \pi^*(\pi_2^*(\omega^{k-q})_{|\Delta}).$$ 
If we identify $\Delta$ with $X$ in the canonical way, then the last current is equal to $R\wedge \pi^*(\omega^{k-q})$ because 
$\pi_2^*(\omega^{k-q})_{|\Delta}$ is identified with $\omega^{k-q}$. Finally, since we assumed that the h-dimension of $R$ is $\geq k-q$, the last current doesn't vanish. It follows that the density between $T\otimes S'$ and $\Delta$ doesn't vanish. This is a contradiction we are looking for.
\endproof

The following corollary will be used in the study of meromorphic maps on surfaces. It  gives a complete characterization of when the density between  two positive closed $(1,1)$-currents is of h-dimension $0$.

\begin{corollary} \label{C:h-surface}
Let $X$ be a compact K\"ahler surface. Let $T$ and $S$ be two positive closed $(1,1)$-currents on $X$. Assume that there is no compact analytic curve $Y$ in $X$ such that both $T$ and $S$ have positive mass on $Y$. Then the density between $T$ and $S$ is of h-dimension $0$. 
\end{corollary}
\proof
Note that by a theorem of Siu, $T$ has positive mass on $Y$ if and only if it has positive Lelong number at each point of $Y$ or equivalently it is 
equal to the sum of a positive closed current and a positive constant times $[Y]$,  see e.g. \cite{Demailly}. Moreover, the set $\{x\in X,\ \nu(T,x)>0\}$ is a finite or countable union of proper analytic subsets of $X$, i.e., of compact analytic curves or points. By hypothesis, $S$ has no mass on the analytic curves in  $\{x\in X,\ \nu(T,x)>0\}$. Since
$S$ cannot has mass at any point of $X$, it has no mass on any countable set. We conclude that $S$ has no mass on  $\{x\in X,\ \nu(T,x)>0\}$. Proposition \ref{P:tensor_prod_hdim} implies the result. 
\endproof
We see that the last result holds even when  the wedge-product $T\wedge S$ is not defined. This somehow illustrates the flexibility of h-dimension and tangent currents in applications. Note also that when both $T$ and $S$ have positive mass on a curve $Y$, it is not difficult to see that the density between $T$ and $S$ is of the maximal h-dimension, i.e., 1.

 \section{Exponential growth rate for isolated periodic points} \label{S:growth_rate}

 In this section, we will give the proofs of Theorems \ref{th_main_1} and \ref{th_main_2}. We need to recall some properties of meromorphic maps on a compact K\"ahler manifold, see also \cite{DinhSibony04,DinhSibony05}. 

Consider a dominant meromorphic map $f:X\to X$ on a compact K\"ahler manifold of dimension $k$ as in the beginning of Introduction. 
We will use the notations already introduced there. 
The current $(f^n)^*(\omega^p)$ is positive closed and of bi-degree $(p,p)$. It is also an $L^1$ form which is  smooth outside the indeterminacy set $I(f^n)$ of $f^n$. 
Recall that the mass of this current 
is comparable with the norm of $(f^n)^*$ on $H^{p,p}(X,\C)$ and we can compute the dynamical degree of order $p$ of $f$ by 
\begin{equation}  \label{E:dp}
d_p(f)=\lim_{n\to\infty} \|(f^n)^*(\omega^p)\|^{1/n}.
\end{equation}
We have the following lemma where we use the natural K\"ahler metric on $X\times X$ associated with the K\"ahler form $\pi_1^*(\omega)+\pi_2^*(\omega)$.

\begin{lemma} \label{L:volume}
Let $\Gamma(f^n)$ denote the closure of the graph of $f^n$ in $X\times X$ and $\|\Gamma(f^n)\|$ the mass of the current of integration on it.
 Then 
$$\lim_{n\to\infty} \|\Gamma(f^n)\|^{1/n} =e^{h_a(f)}=\max_{0\leq p\leq k} d_p(f).$$
\end{lemma}
\proof
We have that
$$\|\Gamma(f^n)\| = \int_{\Gamma(f^n)} \big(\pi_1^*(\omega)+\pi_2^*(\omega)\big)^k=\sum_{p=0}^k {k\choose p} \int_{\Gamma(f^n)} \pi_1^*(\omega^{k-p}) \wedge \pi_2^*(\omega^p).$$
Observe that $\pi_1$ defines a bijective map between a Zariski open set of $\Gamma(f^n)$ and  a Zariski open set of $X$. Moreover,  proper analytic subsets of  $\Gamma(f^n)$ and of $X$ have $2k$-dimensional volume zero. Therefore, by pushing the integrals on $\Gamma(f^n)$ to $X$ by $\pi_1$, we see that the last sum of integrals is equal to
$$\sum_{p=0}^k {k\choose p} \int_X \omega^{k-p} \wedge (f^n)^*(\omega^p).$$
It is now easy to deduce the lemma from \eqref{E:dp}. 
\endproof

Note that by Wirtinger's theorem, the $2k$-dimensional volume of $\Gamma(f^n)$ is  equal to ${1\over k!}\|\Gamma(f^n)\|$.  So in the last lemma we can replace $\|\Gamma(f^n)\|$ by the volume of $\Gamma(f^n)$.  

\bigskip
\noindent {\bf End of the proof of Theorem \ref{th_main_1}.}
Fix an $\epsilon>0.$ We  only need to show that 
$$\lim_{n\to\infty} e^{-n (h_a(f)+\epsilon)} P_n(f)=0.$$
Define $d_n:= e^{n (h_a(f)+\epsilon)}$. 
By Lemma \ref{L:volume}, we have $d_n^{-1} [\Gamma(f^n)]\to 0$.
We apply Proposition
\ref{prop_intersection} for 
$\Gamma(f^n), X\times X, \Delta, 0$ instead of $\Gamma_n, X, V, T$. 
The constant $c$ there is then equal to 0 and the number $\delta_n$ is equal to $P_n(f)$. 
The desired property follows immediately. Note that  Theorem \ref{th_main_1} still  remains true if $f$ is a meromorphic correspondence, see
\cite{DinhSibony08}  for  related results.
\hfill $\square$
  
\begin{remark} \rm
When $X$ is a projective manifold, we can also obtain the result by moving $\Gamma(f^n)$ in a family of (non-effective) cycles with controlled degrees. This approach can be generalized to maps defined over fields different from $\C$ but cannot be used for K\"ahler manifolds, see also \cite{Tuyen}. For K\"ahler manifolds, one can  use the technique of regularization of currents in \cite{DinhSibony04}. However, a mass control for intersections of currents is needed and this method does not yield better results.
\end{remark}
  
We now prove Theorem \ref{th_main_2}. 
From now on,  assume that $X$ is a compact K\"ahler surface, i.e., $k=2$. 
Let $f:X\to X$ be a meromorphic map as in this theorem which is algebraically stable and of minor topological degree. Write for simplicity
$d_1:=d_1(f)$ and  $d_2:=d_2(f).$ We have $d_1>d_2\geq 1$. 
We recall first  some results by Diller-Favre \cite{DillerFavre} and   Diller-Dujardin-Guedj \cite{DillerDujardinGuedj1}.

There are two positive closed $(1,1)$-currents $T_+$ and $T_-$, called {\it dynamical Green currents}, such that 
$f^*(T^+)=d_1 T^+$, $f_*(T^-)=d_1(T^-)$ and $\{T^+\}\smallsmile \{T^-\}=1$. Moreover, the current $T^+$ has no mass on compact analytic curves in $X$. We will need the following proposition.

\begin{proposition} \label{P:pullback-form}
Let $\alpha$ and $\beta$ be two smooth $2$-forms on $X$, not necessarily positive nor closed. 
Define $c_\alpha:=\langle T^-,\alpha\rangle$ and $c'_\beta:=\langle T^+,\beta\rangle$. Then  
$$\lim_{n\to\infty} d_1^{-n} (f^n)^*(\alpha)= c_\alpha T^+ \quad \text{and} \quad \lim_{n\to\infty} d_1^{-n} (f^n)_*(\beta)= c'_\beta T^-.$$
\end{proposition}
\proof
Consider the first limit with $\alpha$ a $(2,0)$-form. Since $T^-$ is of bi-degree $(1,1)$, we have $c_\alpha=0$. 
Let $\phi$ be any smooth test $(0,2)$-form on $X$. We have by the Cauchy-Schwarz inequality
\begin{eqnarray*}
\big|\langle d_1^{-n}  (f^n)^*(\alpha),\phi\rangle \big|&\leq&  d_1^{-n} \Big|\int (f^n)^*(\alpha)\wedge (f^n)^*(\overline\alpha)\Big|^{1/2} \Big|\int \phi\wedge \overline\phi\Big|^{1/2} \\
& = & d_1^{-n} \Big|\int (f^n)^*(\alpha \wedge \overline\alpha)\Big|^{1/2} \Big|\int \phi\wedge \overline\phi\Big|^{1/2},
\end{eqnarray*}
where the integrals are taken on $X$ outside the indeterminacy set of $f^n$. Since $\alpha\wedge\overline\alpha$ is a smooth form of maximal degree, the first integral in the last product is of order at most $O(d_2^n)$. Thus, the property $d_2<d_1$ implies  that $d_1^{-n} \langle (f^n)^*(\alpha),\phi\rangle$ tends to 0 as $n\to\infty$. It follows that $d_1^{-n}  (f^n)^*(\alpha)\to 0$ and the lemma holds in this case 
because $c_\alpha=0$.

Similarly, the lemma holds for $\alpha$ of bi-degree $(0,2)$. We consider now the remaining case where $\alpha$ is of bi-degree $(1,1)$.
This case was already obtained in \cite[Lemma 3.4]{DillerDujardinGuedj1} when $\alpha$ is the product of a smooth function with a K\"ahler form. It is easy to adapt their proof to any smooth $(1,1)$-form. We can also obtain the general case by writing $\alpha$ as a finite linear combination of forms which are products of smooth functions with K\"ahler forms. To see the last point, we can use a partition of unity to reduce the question to the local setting and then use a suitable family of K\"ahler forms to get a local frame for the vector bundle of $(1,1)$ cotangent vectors of $X$.

We will now deduce the second limit in the proposition from the first one. 
Recall that the current $(f^n)_*(\beta)$ is defined as $(\pi_2)_*\big(\pi_1^*(\beta)\wedge [\Gamma(f^n)]\big)$ and is defined by an $L^1$ form which is smooth outside an analytic set. 
For any smooth 2-form $\alpha$, we have
$$\lim_{n\to\infty}\langle d_1^{-n} (f^n)_*(\beta),\alpha\rangle = \lim_{n\to\infty} \langle \beta, d_1^{-n} (f^n)^*(\alpha)\rangle 
=\langle c_\alpha T^+,\beta\rangle=c'_\beta c_\alpha=c'_\beta\langle T^-, \alpha\rangle.$$
This property holds for every smooth 2-form $\alpha$. We deduce that $d_1^{-n} (f^n)_*(\beta)\to c'_\beta T^-$. 
\endproof

\begin{proposition} \label{P:graphs}
Let $\Gamma(f^n)$ denote the closure of the graph of $f^n$ in $X\times X$. Then the sequence of positive closed $(2,2)$-currents
$d_1^{-n}[\Gamma (f^n)]$ converges to  
$T^+\otimes T^- $ as $n\to\infty$.
 \end{proposition}
\proof 
Denote by $(z^1,z^2)$ a general point in $X\times X$ with $z^1,z^2\in X$.
Consider smooth test $(2,2)$-forms  $\Phi$ on $X\times X$ that can be written as
$$\Phi(z^1,z^2) = \beta(z^1)\wedge \alpha(z^2),$$
where $\beta, \alpha$ are  smooth forms on $X$ of bi-degrees respectively $(p,q)$ and $(2-p,2-q)$ with $0\leq p,q\leq 2$.
Since these forms span a dense vector space in the space of test $(2,2)$-forms, we only need to check that 
\begin{equation} \label{E:graphs}
 d_1^{-n} \langle [\Gamma(f^n)],\Phi\rangle  \to \big\langle T^+\otimes T^-,\Phi \big\rangle.
\end{equation}
Note that the left hand side of the last  line is equal to
$$d_1^{-n} \langle [\Gamma(f^n)],\Phi\rangle = d_1^{-n}\int_{\Gamma(f^n)} \pi_1^*(\beta)\wedge \pi_2^*(\alpha) = d_1^{-n}\int_X \beta\wedge (f^n)^*(\alpha)=d_1^{-n}\int_X (f^n)_*(\beta)\wedge \alpha.$$
The right hand side of \eqref{E:graphs} is equal to $\langle T^+,\beta\rangle \langle T^-,\alpha\rangle$. So it vanishes 
except in the first case considered below. 

\medskip
\noindent
{\bf Case 1.}
Assume that $p+q=2$. We have 
$$d_1^{-n} \langle [\Gamma(f^n)],\Phi\rangle =  d_1^{-n}\int_X \beta\wedge (f^n)^*(\alpha)=\langle d_1^{-n} (f^n)^*(\alpha),\beta\rangle.$$
Proposition \ref{P:pullback-form} implies that the last pairing converges to $\langle T^+,\beta\rangle\langle T^-,\alpha\rangle$. So \eqref{E:graphs} holds in this case.

\medskip
\noindent
{\bf Case 2.} Assume that $p=q=0$. So $\beta$ is a function and $\alpha$ is a form of maximal degree. 
In particular, $\alpha$ defines a measure and therefore the mass of $(f^n)^*(\alpha)$ is of order at most $O(d_2^n)$. We then deduce that 
$$d_1^{-n} \langle [\Gamma(f^n)],\Phi\rangle =  \langle d_1^{-n} (f^n)^*(\alpha),\beta\rangle$$
is of order at most $O(d_1^{-n} d_2^n)$. The assumption that $d_2<d_1$ implies that $\langle d_1^{-n} [\Gamma(f^n)],\Phi\rangle\to 0$ and we get \eqref{E:graphs} in this case.

\medskip
\noindent
{\bf Case 3.}  Assume that $p=q=2$. So $\alpha$ is a function and $\beta$ is a form of maximal degree.  It follows that $(f^n)^*(\alpha)$ is a bounded function and we easily deduce that $\langle d_1^{-n} [\Gamma(f^n)],\Phi\rangle = O(d_1^{-n})$ and  get \eqref{E:graphs} in this case.

\medskip
\noindent
{\bf Case 4.} Assume that $p+q=1$. We only consider the case $p=1$ and $q=0$ because the other case with $p=0$ and $q=1$ can be treated in the same way. So $\beta$ is a $(1,0)$-form and $\alpha$ is a $(1,2)$-form. We can assume that $\alpha=\gamma\wedge \theta$, where $\gamma$ is a $(0,1)$-form and $\theta$ is a K\"ahler form. Indeed, $\alpha$ can be written as a finite combination of such forms. 
By the Cauchy-Schwarz inequality, we have 
$$d_1^{-n} |\langle [\Gamma(f^n)],\Phi\rangle|\leq d_1^{-n}  \big |\big\langle [\Gamma(f^n)],\beta(z^1)\wedge \overline{\beta(z^1)}\wedge \theta(z^2)\big\rangle \big|^{1/2}  \big|\big\langle [\Gamma(f^n)],\gamma(z^2)\wedge \overline{\gamma(z^2)}\wedge \theta(z^2)\big\rangle\big|^{1/2}.$$
We can apply the estimates in Cases 1 and 2 in order to bound the factors in the last product. 
We obtain that this product is of order $O(d_1^{-n} d_1^{n/2}d_2^{n/2})$. It follows that 
$d_1^{-n} \langle [\Gamma(f^n)],\Phi\rangle\to 0$ and we obtain \eqref{E:graphs} for this case.

\medskip
\noindent
{\bf Case 5.} Assume in this remaining case that $p+q=3$. We can also suppose that $p=2$, $q=1$ and $\beta$ equal to the wedge-product of a $(1,0)$-form with a K\"ahler form. Then, this case can be treated exactly as in the last case using the estimates from Cases 1 and 3. 
\endproof

\smallskip

\noindent{\bf  End of the  proof of Theorem \ref{th_main_2}.}
By Proposition \ref{P:graphs}, we have $d_1^{-n} [\Gamma(f^n)]\to T^+\otimes T^-$. 
Since $T^+$ has no mass on compact analytic curves in $X$, we can apply Corollary \ref{C:h-surface} to $T^+, T^-$ instead of $T$ and $S$.
We obtain that the tangential h-dimension between $T^+\otimes T^-$ along $\Delta$ is 0. 
Therefore, we can apply Proposition \ref{prop_intersection} to $X\times X, \Delta, d_1^n, \Gamma(f^n), T^+\otimes T^-$ instead of $X,V,d_n, \Gamma_n$ and $T$. The constant $c$ in this proposition is equal to 
$$\{T^+\otimes T^-\}\smallsmile \{\Delta\}= \{T^+\}\smallsmile \{T^-\}=1.$$
The constant $\delta_n$ there is equal to $P_n(f)$. Theorem \ref{th_main_2} follows immediately. 
\hfill $\square$

\medskip

Recall that the Lefschetz number $L(f^n)$ is the intersection number $\{\Gamma(f^n)\}\smallsmile \{\Delta\}$. We have the following useful lemma.

\begin{lemma} \label{L:Lefschetz}
Under the hypotheses of Theorem \ref{th_main_2} and with the above notations, we have 
$$L(f^n)=d_1^n+o(d_1^n)$$
as $n\to\infty$. 
\end{lemma}
\proof
Let $\alpha_{p,q,i}$ be smooth closed $(p,q)$-forms on $X$ with $1\leq i\leq \dim H^{p,q}(X,\C)$ such that the classes $\{\alpha_{p,q,i}\}$ give a basis of the Hodge cohomology group $H^{p,q}(X,\C)$. We also choose these forms so that the basis for $H^{p,q}(X,\C)$ is dual to the basis for $H^{2-p,2-q}(X,\C)$, i.e., $\{\alpha_{p,q,i}\}\smallsmile \{\alpha_{2-p,2-q,j}\}$ is equal to 0 when $i\not=j$ and 1 if $i=j$. So the class of $\Delta$ in $H^{2,2}(X\times X,\C)$ is equal to the class of 
$$\sum_{p,q,i} \alpha_{p,q,i}(z^1)\wedge \alpha_{2-p,2-q,i}(z^2).$$
We then deduce that
$$L(f^n)=\sum_{p,q,i} \int_{\Gamma(f^n)} \alpha_{p,q,i}(z^1)\wedge \alpha_{2-p,2-q}(z^2)=\sum_{p,q,i} \int_X \alpha_{p,q,i}\wedge (f^n)^*(\alpha_{2-p,2-q,i}).$$
We have seen in the estimates from the proofs of Propositions \ref{P:pullback-form} and \ref{P:graphs} that the last integral is of order $o(d_1^n)$ when $(p,q)\not=(1,1)$. 

The rest in the sum of integrals (those corresponding to $p=q=1$) can be interpreted as the trace of $(f^n)^*$ on $H^{1,1}(X,\C)$. Recall that since $f$ is algebraically stable we have $(f^n)^*=(f^*)^n$ on $H^{1,1}(X,\C)$. Moreover, it is known from  \cite{DillerFavre}  that $d_1$ is a simple eigenvalue of $f^*$ acting on $H^{1,1}(X,\C)$ and the other eigenvalues have moduli smaller than $d_1$. It is now clear that the above trace of $(f^n)^*$ is equal to $d_1^n+o(d_1^n)$. The lemma follows.  Note that a more precise estimate in terms of $d_1, d_2$ and the eigenvalues of $f^*$ on $H^{1,1}(X,\C)$ can be obtained. 
\endproof

\section{Bi-meromorphic maps on compact K\"ahler surfaces} \label{S:bimero_maps} 

In this section, we will show that the theory of densities of currents may be combined with other arguments to get a lower estimate for the number of isolated periodic points. Namely, we will give here the proof of Theorem \ref{th_main_3}. From now on, we assume that $f:X\to X$ is an algebraically stable  bi-meromorphic map on a compact K\"ahler surface as in that theorem. 

We  start with  a  short digression on Iwasaki-Uehara's version of Saito's fixed point formula in our setting. 
Denote by $I(f^n)$ and $I(f^{-n})$ the indeterminacy sets of $f^n$ and $f^{-n}$ respectively. These sets are finite. 
Recall that since $f$ is algebraically stable, these two sets are disjoint, see e.g. \cite[Proposition 1.13]{DillerFavre} and \cite[Theorem 2.4]{DillerJacksonSommese}. So if $x$ is a periodic point of period $n$ of $f$, then either $f^n$ is holomorphic near $x$ and fixes the point $x$ or the similar property holds for $f^{-n}$. 
Saito's local index function  is the  function $\nu_\bullet(f^n):\Per_n(f)\to\N$ defined in  \cite[Definition 3.5]{IwasakiUehara}. This function is $0$ except at a finite number of points. 
We only recall this notion in the case of isolated periodic points, since this is sufficient for the purpose here. 
 
Let $x$ be an isolated fixed point of $f$ and for simplicity assume that $f$ is holomorphic near $x$.  
Fix a local coordinate $z=(z_1,z_2)$ of $X$  centered at $x$, i.e., $(z_1,z_2)=(0,0)$ at $x$. In these coordinates, we can write
$$f(z)=(z_1+ h_1(z_1,z_2), z_2+h_2(z_1,z_2)),$$
where $h_1$ and $h_2$ are holomorphic functions on an open  neighborhood of $0$ in $\C^2$, vanishing at 0. 
  Let $\Oc_0$  denote  the ring of germs of holomorphic functions  at $0\in\C^2.$
  Let $\mathcal I\subset \Oc_0$ be the  ideal generated by $h_1$ and $h_2.$
  Since $0$ is an isolated fixed point of $f,$  the germ of analytic set defined by $ \mathcal I $ coincides with the single point $ \{0\}.$
  In other words,  $ \mathcal I $ is of complete intersection. By  \cite[Definition 3.5]{IwasakiUehara}, 
$$\nu_{x}(f):=\dim_\C  \Oc_0/\mathcal I. $$
The following lemma shows that Saito's index extends the notion of multiplicity for isolated periodic points.

\begin{lemma}
If $x$ is an isolated fixed point of $f$ as above, then $\nu_x(f)$ is the multiplicity of the intersection $\Gamma(f)\cap \Delta$ at the point $(x,x)$. 
\end{lemma}
\proof
Let $d$ denote the  multiplicity of the intersection $\Gamma(f)\cap \Delta$ at the point $(x,x)$.
 Consider the germ of proper holomorphic map $h:=(h_1,h_2):\ (\C^2,0) \to (\C^2,0)$.
 This is a ramified  covering of degree $d$ from a neighborhood $U$ of $0$ in $\C^2$ to $h(U)$.
 The degree $d$  represents the multiplicity of $\mathcal I$ at $0,$ see e.g.  \cite{Mumford}.
Since $\mathcal I$ is of complete intersection, this multiplicity is  equal to the index $\nu_x(f)$ defined above,  see e.g.   \cite{Serre}.
The lemma follows.
\endproof

Recall now the index of a fixed  curve, i.e.,  a curve of fixed points. The notion can be easily extended to curves of periodic points.  Let $C$ be an  irreducible analytic  curve  in $X$. We say that $C$  is  {\it a fixed curve} of $f$  if  $f(x)=x$ for all $x\in C\setminus I(f)$ or equivalently 
$C$ is contained in $\Per_1(f)$. 
Let $X_1(f)$ be the set of all fixed curves of $f.$  Following \cite[Definition 3.5]{IwasakiUehara}, we can define the index function
 $\nu_\bullet (f):\   X_1(f)\to\N$  as follows, see also \cite[Lemma 6.1]{IwasakiUehara}. Let $C$ be a curve in $X_1(f).$ Let $x$ be a regular point of $C$ which is not a point of indeterminacy for $f$.
 We can  find  local coordinates $z=(z_1,z_2)$ of $X$ centered at $x$  such that  in 
these coordinates
$$
f(z)=(z_1+ z_1^p h_1, z_2+z_1^qh_2),
$$
where $p,q\in\N_+,$ $h_1,h_2$ are holomorphic on a neighborhood of $0,$
  and $h_1(0,z_2), h_2(0,z_2)$  do not vanish identically.
Define
$$\nu_C(f):= \min(p,q).$$
This notion is independent of the choices of coordinates. 
Observe that except for a finite number of points $x$, we have $h_1(0)\not=0$ and $h_2(0)\not=0$, and therefore, reducing the neighborhood of $x$ allows us to assume that $h_1$ and $h_2$ are nowhere vanishing.

The following lemma relates $\nu_C(f)$ to the density of the current $[\Gamma(f)]$ along $\Delta$, see Section \ref{S:currents} for notations. Recall that we identify $\Delta$ with $X$ in a canonical way.

\begin{lemma}\label{L:comparison}
The  tangential h-dimension of $[\Gamma(f)]$ along $\Delta$  is  at most $1$ and we have
$$\kappa_1^\Delta([\Gamma(f)])= \sum_{C\in X_1(f)}\nu_C(f)\{C\}.$$ 
 \end{lemma}
 \proof 
For simplicity, write $\Gamma:=\Gamma(f)$, $X_1:=X_1(f)$, $\nu_C:=\nu_C(f)$ and $\kappa_j:=\kappa^\Delta_j([\Gamma(f)])$.
We will apply the theory of densities of currents to $X\times X, [\Gamma(f)], \Delta$ instead of $X,T, V$, see Section \ref{S:currents} for details. We will use notations similar to those in that section. Denote by $\overline E$ the natural compactification of the normal vector bundle $E$ of  $\Delta$ in $X\times X$ and 
$\pi:\overline E\to \Delta$ the canonical projection. Consider an arbitrary tangent current $R$ of $[\Gamma]$ along $\Delta$. The map $A_\lambda:\overline E\to \overline E$ is the multiplication by $\lambda$ on the fibers of $\overline E$. The first assertion in the lemma is automatically true and we only have to prove the second one.

Recall that if $-h$ is the tautological $(1,1)$-class of $\overline E$, we have the following unique decomposition 
$\{R\}=\pi^*(\kappa_1)\smallsmile h + \pi^*(\kappa_0)$, see \eqref{e:total_class}.  
We have to show that 
$$\kappa_1= \sum_{C\in X_1}\nu_C\{C\}.$$
Let $\alpha$ be an arbitrary  closed $(1,1)$-form on $\Delta$. By Poincar\'e's duality, it is enough to check that 
 \begin{equation} \label{e:kappa_alpha}
 \kappa_1\smallsmile \{\alpha\}= \sum_{C\in X_1}\nu_C\{C\}\smallsmile\{\alpha\}.
 \end{equation}
We can assume that $\alpha$ is a K\"ahler form because the classes of such forms spans $H^{1,1}(X,\C)$. 
So the right hand side of the last identity is the mass of the positive measure $\sum \nu_C [C]\wedge \alpha$. Denote this mass by $m$. 

Observe that $\{R\wedge \pi^*(\alpha)\}=\pi^*(\kappa_1\smallsmile\{\alpha\})\smallsmile h$. As we have seen in the discussion before \eqref{e:total_class}, this decomposition is unique. Thus, identity \eqref{e:kappa_alpha} is equivalent to
\begin{equation} \label{e:kappa_1}
\{R\wedge \pi^*(\alpha)\}=mH,
\end{equation}
where $H$ denotes the restriction of $h$ to a fiber of $\overline E$. Note that $H$ is the class of a hyperplane of a fiber of $\overline E$.

Consider a point $(x,x)\in \Delta$ with $x\in X$. Let $z=(z_1,z_2)$ denote a local coordinate system in $X$ centered at $x$, which identifies a neighborhood of $x$ with a neighborhood $V_0$ of $0$ in $\C^2$. It also induces naturally a local coordinate system $(z,z')$ with $z'=(z_1',z_2')$ of $X\times X$ centered at $(x,x)$. The diagonal $\Delta$ is given by $z=z'$ in these coordinates. Define $w=(w_1,w_2)=z'-z$. So $(z,w)$ is a new local coordinate system centered at $(x,x)$ and $\Delta$ is given by $w=0$. We can identify a neighborhood of $(x,x)$ with some neighborhood $U$ of $V_0\times \{0\}$ in $V_0\times \C^2$. Here, $\Delta$ is identified with $V_0\times\{0\}$, the normal bundle $E$ with $V_0\times \C^2$, and the projection $\pi$ with the natural projection from $V_0\times\C^2$ to $V_0$. The map $A_\lambda$ is then given by $A_\lambda(z,w)=(z,\lambda w)$. We also choose the admissible map $\tau$ to be the identity map as the standard local setting described in Section \ref{S:currents}.

We see in this picture that $(A_\lambda)_*[\Gamma]$ is given by the analytic subset $A_\lambda(\Gamma)$ of pure dimension 2 in $A_\lambda(U)$. The last set is an open set converging to $V_0\times \C^2$ as $\lambda\to\infty$. Given any compact set in $V_0\times \C^2$, the theory of densities insures that the volume of $A_\lambda(\Gamma)$ (or equivalently the mass of the associated current) in $K$ is bounded when $\lambda\to\infty$. This allows us to extract a converging subsequence from the family of currents $(A_\lambda)_*[\Gamma]$ and get a tangent current $R$, see Section \ref{S:currents}. We deduce that $R$ is given by a finite linear combination $\sum m_i\Gamma_i$  of irreducible sub-varieties $\Gamma_i$ of dimension 2 in $E$, where $m_i$ are positive integers. We refer to King \cite{King} for details on currents defined by positive chains of varieties. 
The varieties  $\Gamma_i$ are invariant by $A_\lambda$ for every $\lambda$ because $R$ satisfies this property. 
 
If $x$ is not a fixed point, we can choose $U$ small enough so that it has no intersection with $\Gamma$. In this case, we see that $R=0$ in $V_0\times \C^2$. In contrast, if $x$ is a fixed point, then $A_\lambda(\Gamma)$ contains $(x,x)$ and therefore some $\Gamma_i$ contains this point. We conclude that the union of 
$\Gamma_i\cap\Delta$ is exactly the set of points $(x,x)$ with $x\in\Per_1(f)$. 
Consider first an index $i$ such that $\Gamma_i\cap\Delta$ is finite. Since $\Gamma_i$ is irreducible, of dimension 2, and invariant by $A_\lambda$, we deduce that $\Gamma_i$ is a fiber of $E$. Therefore, we get $[\Gamma_i]\wedge \pi^*(\alpha)=0$. 
So such a $\Gamma_i$ does not contribute to the left hand side of  identity \eqref{e:kappa_1}.
We then deduce that $R\wedge \pi^*(\alpha)$ has no mass on each fiber of $\pi$. It remains now to analyze the situation near a  point $(x,x)$ of $\Delta$ where $x$ is a generic point on a curve $C$ in $X_1$.

Recall that tangent currents do not depend on the choice of coordinates. So we can choose $z$ as in the discussion just before Lemma 
\ref{L:comparison}. With this choice, the variety $A_\lambda(\Gamma)$ is given by
$$w_1=\lambda z_1^ph_1(z_1,z_2) \quad \text{and} \quad w_2=\lambda z_1^qh_2(z_1,z_2).$$
Consider first the case where $p=q$. The last variety is a ramified covering of degree $p$ over a neighborhood of 0 in the plane $(z_2,w_1)$. 
Taking $\lambda\to\infty$, we see that the current $(A_\lambda)_*([\Gamma])$ converges to $p$ times the current of integration on the variety
$$z_1=0\quad \text{and} \quad w_2={h_2(0,z_2)\over h_1(0,z_2)} w_1.$$
When $p<q$, the limit is $p$ times the current of integration on $z_1=w_2=0$ and when $p>q$ the limit is $q$ times the current of integration on $z_1=w_1=0$. In any case, we see that $R\wedge \pi^*(\alpha)$ restricted to $\pi^{-1}(V_0)$ can be decomposed into currents of integration on hyperplanes of the fibers of $\pi$. Moreover, its cohomology class is equal to the mass of $\nu_C\cdot [C]\wedge\alpha$ in $V_0$ times $H$. The desired identity \eqref{e:kappa_1} follows.
 \endproof

Recall that
the set $X_1(f)$  of all fixed curves can be divided into
 two disjoint families, namely,
 into what Saito called {\it the curves of type I} and {\it type II}:
$$X_1( f ) = X_\rI ( f ) \sqcup X_\rII ( f ). $$
With the notation as above, the curve $C$ is of type I when $p\leq q$ and of type II when $p>q$. One can check that the definition does not depend on the choice of coordinates.
Roughly speaking, $C$ is of type II if near a generic point of $C$, the map $f$ moves slightly a point 
to another along a tangential direction to $C$. 
We refer the reader    to \cite[Def. 3.7 and Lemma 6.1]{IwasakiUehara} for more details. 

Now  we are in the position to state  Iwasaki-Uehara's version of  Saito's fixed point formula.
This is a key argument in the proof of Theorem \ref{th_main_3}. 
Note that  these authors state their result for projective surfaces but their proof also works for K\"ahler surfaces. Of course, this result also applies to $f^n$ with $n\in\Z$. 
 
\begin{theorem}\label{T:Saito}  {\rm  (Saito  \cite{Saito}, Iwasaki-Uehara \cite{IwasakiUehara})}
Under the  above hypotheses and  notations,  the Lefschetz number  $L(f)$ can be  expressed as 
$$
L(f)=\sum _{x\in \Per_1(f)}\nu _x(f)+\sum _{C\in X_\rI(f)}\chi _C\cdot \nu _C(f)+\sum _{C\in X_\rII(f)}\tau_C\cdot\nu _C(f),
$$
where $\chi_C$ is the  Euler characteristic of the  normalization of $C$ and $\tau_C$ is the self-intersection number of $C.$
 \end{theorem}

\noindent{\bf End of the proof of Theorem \ref{th_main_3}.} 
By Theorem \ref{T:Saito}, applied to $f^n$, we have
\begin{equation}\label{e:Saito_n}
L(f^n)=\sum _{x\in \Per_n(f)}\nu _x(f^n)+\sum _{C\in X_\rI(f^n)}\chi _C\cdot \nu _C(f^n)+\sum _{C\in X_\rII(f^n)}\tau_C\cdot
\nu _C(f^n).
\end{equation}
We need to show that the first term in the right hand side is equal to $d_1(f)^n+o(d_1(f)^n)$. By Lemma \ref{L:Lefschetz}, we only need to check that the last two terms in the previous identity is equal to $o(d_1(f)^n)$. 

Consider the first term with $C\in X_\rI(f^n)$. A  result of Diller-Jackson-Sommese \cite[Theorem 3.6] {DillerJacksonSommese} says that
the genus of such a curve  $C$ is  zero or one, that is, $\chi _C$ is $2$ or $0.$
Therefore, we have
\begin{equation} \label{e:sum_1}
0\leq \sum _{C\in X_\rI(f^n)}\chi _C\cdot \nu _C(f^n)\leq 2 \sum _{C\in X_1(f^n)}  \nu _C(f^n).
\end{equation}
For simplicity, denote by $\Gamma_n$  the closure of the graph of $f^n$ in $X\times X$.
By Lemma \ref{L:comparison}, the  tangential h-dimension of $[\Gamma_n]$ along $\Delta$  is at most equal to $1$ and 
$$\kappa_1^\Delta([\Gamma_n])= \sum_{C\in X_1(f^n)}\nu_C(f^n)\{C\}.$$ 
Moreover,  we know by  Proposition \ref{P:graphs} 
that  $d_1(f)^{-n}[\Gamma_n]\to T^+\otimes T^-$ and by Corollary \ref{C:h-surface} that
$T^+\otimes T^-$  has the tangential h-dimension 0 along $\Delta.$
Consequently, applying Proposition \ref{P:semicontinuity} yields
$$d_1(f)^{-n} \kappa_1^\Delta([\Gamma_n])\to 0 \quad \text{or \ equivalently} \quad 
\sum_{C\in X_1(f^n)}\nu_C(f^n)\{C\} =o(d_1(f)^n).$$
Thus,
\begin{equation} \label{e:sum_nu}
\sum_{C\in X_1(f^n)}\nu_C(f^n)\{C\}\smallsmile\{\omega\}=o(d_1(f)^n),
\end{equation}
where $\omega$ is any fixed K\"ahler form on $X.$
Note that $\{C\}\smallsmile\{\omega\}$ is 2 times the area of $C$ with respect to the K\"ahler metric $\omega$. Therefore, it is bounded below by a positive constant independent of $C$. This together with 
\eqref{e:sum_1} implies that the second term in the right hand side of \eqref{e:Saito_n} is equal to $o(d_1(f)^n)$.

It remains to prove a similar bound for the last term in  \eqref{e:Saito_n}. For this purpose, we have two different approaches. 

\medskip\noindent
{\bf First approach.}  We know  by  Iwasaki-Uehara \cite[Theorem 2.4]{IwasakiUehara} that
 $\mathcal X:=\cup_{n=1}^\infty X_\rII(f^n)$ is a finite set. For every $C\in \mathcal X,$ there is a smallest integer $k_C$  such that
 $C\in  X_\rII(f^{k_C}).$
 Moreover,  if  $C\in X_\rII(f^n),$    \cite[Theorem 2.1]{IwasakiUehara} tells us that 
$\nu_C(f^{kn})=\nu_C(f^n)$ for all $k\geq 1.$
Consequently, for every $n\geq 1$ and $C\in X_\rII(f^n),$
we have
$$\nu_C(f^n)=\nu_C(f^{nk_C})=\nu_C(f^{k_C}).$$
Thus,
$$\Big|\sum _{C\in X_\rII(f^n)}\tau_C\cdot
\nu _C(f^n)\Big|\leq   \sum _{C\in \mathcal X}|\tau_C|\cdot
\nu _C(f^{k_C}).$$
The last sum is a constant independent of $n$. The desired estimate follows.

\medskip\noindent
{\bf Second approach.}  
By Diller-Jackson-Sommese \cite[Theorem 3.6] {DillerJacksonSommese},
  the arithmetic genus $g(C)$ of $C\in X_1(f^n)$  is either zero or one. By the genus formula (e.g. Section 11, Chapter 2 in \cite{BarthPetersVen}),  $g(C)$ and the self-intersection $\tau_C$ of a curve $C$ in the surface $X$ are related via
$$\tau_C=2g(C)-2-K_X\smallsmile \{C\},$$
where  $K_X$ is the canonical
class of $X.$ Consequently,  there is a constant $M>0$  independent of $C$ such that
$$
|\tau_C|\leq M\cdot  \{C\} \smallsmile\{\omega\}.
$$  
Therefore,
\begin{eqnarray*}
\Big|\sum _{C\in X_\rII(f^n)}\tau_C\cdot
\nu _C(f^n)\Big|
&\leq& M\sum_{C\in X_1(f^n)}\nu_C(f^n)\cdot\{C\}\smallsmile\{\omega\}.
\end{eqnarray*}
We then conclude the proof    using   estimate \eqref{e:sum_nu}. 
\hfill $\square$

\end{document}